\documentclass[a4paper,12pt]{amsart}
\usepackage{graphicx}
\newtheorem{theorem}{Theorem}[section]
\newtheorem{proposition}[theorem]{Proposition}
\newtheorem{lemma}[theorem]{Lemma}
\newtheorem{claim}[theorem]{Claim}

\theoremstyle{definition}

\theoremstyle{remark}

\numberwithin{equation}{section}

\begin{document}

\title[Dehn surgeries on knots which yield lens spaces]
{Dehn surgeries on knots which yield lens spaces and genera of knots}

\author{Hiroshi Goda}
\address{The Graduate School of Science and Technology,
Kobe University, Rokkodai 1-1, Nada, Kobe 657-8501, Japan}
\email{goda@math.kobe-u.ac.jp}
\thanks{The first author was supported in part by Hyogo Science 
and Technology Association.}

\author{Masakazu Teragaito}
\address{Department of Mathematics and Mathematics Education,
Hiroshima University, Kagamiyama 1-1-1, Higashi-hiroshima 739-8524, Japan}
\email{mteraga@sed.hiroshima-u.ac.jp}
\thanks{The second author was supported in part by Grant-in-Aid for
Encouragement of Young Scientists 09740059, The Ministry of Education,
Science and Culture.}




\begin{abstract}
Let $K$ be a hyperbolic knot in the 3-sphere.
If $r$-surgery on $K$ yields a lens space, then we show
that the order of the fundamental group of the lens space is
at most $12g-7$, where $g$ is the genus of $K$.
If we specialize to genus one case,
it will be proved that no lens space can be obtained from genus one, hyperbolic knots
by Dehn surgery.
Therefore, together with known facts,
we have that a genus one knot $K$ admits Dehn surgery yielding
a lens space if and only if $K$ is the trefoil.
\end{abstract}

\maketitle

\section{Introduction}

It is well known that every closed orientable $3$-manifold can be obtained
by Dehn surgery on a link in the $3$-sphere $S^3$ \cite{Li,Wal}.
When one considers Dehn surgery on knots, it is natural to think that
there are some restrictions on the resulting manifolds after Dehn surgery,
aside from obvious ones, such as the weight of the fundamental group, or homology groups.
The fact that a knot is determined by its complement \cite{GL1} can
be expressed that a non-trivial surgery on a non-trivial knot never yield
$S^3$.  Similarly, the Property R conjecture solved in \cite{Ga} means that
$S^2\times S^1$ cannot be obtained by surgery on non-trivial knots.
These two results suggest that it is hard to obtain a $3$-manifold
with a relatively simple structure in view of Heegaard genera
by Dehn surgery on knots. 

A lens space $L(m,n)$ is the manifold of Heegaard genus one, and
it can be obtained by $m/n$-surgery on a trivial knot.
There are many studies on the problem of what kind of a knot in
$S^3$ admits Dehn surgery
yielding a lens space.
For torus knots and satellite knots, the question of when Dehn surgery 
on such a knot yields a
lens space is completely solved \cite{BL,Mo,Wa,Wu}.
More precisely, any torus knot admits an infinitely many surgeries yielding
lens spaces, and only the $(2pq\pm 1,2)$-cable knot of the $(p,q)$-torus knot
admits such surgery.

It is also known that there are many examples of hyperbolic knots which
admit Dehn surgery yielding lens spaces.
For example, Fintushel-Stern \cite{FS} have shown that
$18$- and $19$-surgeries on the $(-2,3,7)$-pretzel
knot give lens spaces $L(18,5)$ and $L(19,7)$, respectively.
We note that the $(-2,3,7)$-pretzel knot has genus 5.
As far as we know, this is the minimum among the genera of such 
hyperbolic knots.  
All known examples can be expressed as closed positive (or negative) braids, and
therefore they are fibered \cite{St}, and it is easy to calculate their genera,
since Seifert's algorithm gives fiber surfaces, that is, minimal genus
Seifert surfaces for such knots.
It is conjectured that if a hyperbolic knot
admits Dehn surgery yielding a lens space then the knot is fibered.

On the other hand, Berge \cite{Be} gave a list
in which the knot admits Dehn surgery yielding a lens space.
It is expected \cite{Go0} that this will give a complete list of knots with
Dehn surgery yielding lens spaces, but there seems to be no essential progress
in this direction yet.

In the opposite direction, several families of hyperbolic knots are known to
admit no surgery yielding lens spaces: 
$2$-bridge knots \cite{Ta}, alternating knots \cite{DR}, some Montesinos knots \cite{BZ}.

In this paper we focus on the genera of knots as a new standpoint, and show that
there is a constraint on the order of the fundamental group of
the resulting lens space obtained by surgery on a hyperbolic knot.

Let $K$ be a hyperbolic knot in $S^3$. 
The exterior of $K$, denoted by $E(K)$, is the complement of an open
tubular neighborhood of $K$.  Let $r$ be a \textit{slope\/} on $\partial E(K)$,
that is, the isotopy class of
an essential simple closed curve in $\partial E(K)$, and let
$K(r)$ be the closed 3-manifold obtained by $r$-Dehn surgery on $K$.
That is, $K(r)=E(K)\cup V_r$, where $V_r$ is a solid torus
attached to $\partial E(K)$
along their boundaries in such a way that $r$ bounds a meridian disk in $V_r$.
Slopes on $\partial E(K)$ are parameterized 
as $m/n\in \mathbb{Q}\cup\{1/0\}$ in the usual way \cite{Ro}.

If $K(r)$ is a lens space, then $r$ is an integer
by the Cyclic Surgery Theorem \cite{CGLS}.
In particular, $\pi_1K(r)$ has the order $|r|$.

\begin{theorem}
Let $K$ be a hyperbolic knot in $S^3$.
If $K(r)$ is a lens space,
then $|r|\le 12g-7$, where $g$ denotes the genus of $K$.
\end{theorem}

For genus one case, 
we have the complete answer.

\begin{theorem} No Dehn surgery on a genus one, hyperbolic knot in
$S^3$ gives a lens space.
\end{theorem}

Combining this with known facts, we can completely determine Dehn surgeries
on genus one knots which yield lens spaces.

\begin{theorem}
A genus one knot $K$ in $S^3$ admits
Dehn surgery yielding a lens space
if and only if $K$ is the $(\pm 3,2)$-torus knot
and the surgery slope is $(\pm 6n+\varepsilon)/n$ for $n\ne 0$, $\varepsilon=\pm 1$.
\end{theorem}

As earlier results, we have proved that no Dehn surgery on a genus one
knot gives $L(2,1)$ \cite{Te1} (see also \cite{Do}) and $L(4k,2k\pm 1)$ for
$k\ge 1$ \cite{Te2}. 
It was also known that if a genus one knot has a non-trivial Alexander polynomial,
then the knot has no cyclic surgery of even order \cite[Corollary 2]{Ne}. 
Recently, \cite{HS} showed that the lens space
$L(2p,1)$ cannot be obtained by surgery on a strongly invertible knot.

To prove Theorems 1.1 and 1.2, we will analyze the graphs of the intersection of the punctured
surfaces in a knot exterior coming from
a Heegaard torus of a lens space and a minimal genus Seifert surface for the knot.
By virtue of the use of a Seifert surface,
instead of a level sphere in a thin position of the knot,
the graphs can include
the information on the order of
the fundamental group of the resulting lens space after Dehn surgery.
In Section 2, it will be found out that there are some constraints on
Scharlemann cycles.
The proof of Theorem 1.1 is divided into two cases according to the number $t$ of points
of intersection between the Heegaard torus and the core of the attached solid torus.
In Section 3, the case that $t\ge 4$ is dealt with, and the special case that
$t=2$ is discussed in Section 4 and the proof of Theorem 1.1 is completed.
Finally in Sections 5 and 6, we specialize to the case that $K$ has genus one,
and prove Theorems 1.2 and 1.3.

\section{Preliminaries}

Throughout this paper, $K$ will be assumed to be a hyperbolic knot in $S^3$.
For a slope $r$, suppose that $K(r)=E(K)\cup V_r$ is a lens space.
Since $K$ is not a torus knot, the Cyclic Surgery Theorem \cite[Corollary 1]{CGLS}
implies that the slope $r$ must be integral.
We may assume that $r>1$.
Thus $\pi_1K(r)$ has the order $r$.
For simplicity, we denote $V_r$ by $V$.
Let $K^*$ be the core of $V$.

Let $\widehat{T}$ be a Heegaard torus in $K(r)$.
Then $K(r)=U\cup W$, where $U$ and $W$ are solid tori.
We can assume that $\widehat{T}$ meets $K^*$ transversely in
$t$ points, and that $\widehat{T}\cap V$ consists of
$t$ mutually disjoint meridian disks of $V$.
Then $T=\widehat{T}\cap E(K)$ is a punctured torus with $t$ boundary components,
each having slope $r$ on $\partial E(K)$.

Let $S\subset E(K)$ be a minimal genus Seifert surface of $K$.
Then $S$ is incompressible and boundary-incompressible in $E(K)$.

By an isotopy of $S$, we may assume that
$S$ and $T$ intersect transversely, and $\partial S$ meets each
component of $\partial T$ in exactly $r$ points.
We choose $\widehat{T}$ so that the next condition $(*)$
is satisfied :

\medskip
$(*)$ $\widehat{T}\cap K^*\ne \emptyset$, and
each arc component of $S\cap T$ is essential in $S$ and in $T$.
\medskip

This can be achieved if $K^*$ is put in thin position with respect to
$\widehat{T}$ \cite{Ga,Go1}.
(Note that if $K^*$ can be isotoped to lie on $\widehat{T}$, then
$K$ would be a torus knot.)
Furthermore, we may assume that $\widehat{T}$ is chosen so that
$t$ is minimal over all Heegaard tori in $K(r)$ satisfying $(*)$.
This minimality of $\widehat{T}$ will be crucial in this paper.

Since $S$ is incompressible in $E(K)$ and $E(K)$ is irreducible, it can be assumed that
no circle component of $S\cap T$ bounds a disk in $T$.
But it does not hold for $S$ in general.
We further assume that the number of loop components of $S\cap T$ is
minimal up to an isotopy of $S$.

The arc components of $S\cap T$ define graphs $G_S$ in $\widehat{S}$ and $G_T$
in $\widehat{T}$ as follows \cite{CGLS,GL1}, where $\widehat{S}$ is the closed
surface obtained by capping $\partial S$ off by a disk.
Let $G_S$ be the graph in $\widehat{S}$ obtained by
taking as the (fat) vertex the disk $\widehat{S}-\mathrm{Int}S$
and as edges the arc components of $S\cap T$ in $\widehat{S}$.
Similarly, $G_T$ is the graph in $\widehat{T}$ whose vertices are
the disks $\widehat{T}-\mathrm{Int}T$
and whose edges are the arc components of $S\cap T$ in $\widehat T$.
Number the components of $\partial T$, $1,2,\ldots,t$ in sequence
along $\partial E(K)$.  Let $\partial_iT$ denote the component of $\partial T$ with label $i$.
This induces a numbering of the vertices of $G_T$.
Let $u_i$ be the vertex of $G_T$ with the label $i$ for $i=1,2,\dots,t$.
Let $H_{x,x+1}$ is the part of $V$ between
consecutive fat vertices $u_x$ and $u_{x+1}$ of $G_T$.
When $t=2$, $V$ is considered to be the union $H_{1,2}\cup H_{2,1}$.
Each endpoint of an edge in $G_S$ at the unique vertex $v$
has a label, namely the label of the corresponding component of $\partial T$.
Thus the labels $1,2,\ldots,t$ appear in order around $v$ repeated $r$ times.

The graphs $G_S$ and $G_T$ satisfy the parity rule \cite{CGLS}
which can be expressed as the following :
the labels at the endpoints of an edge of $G_S$ have distinct parities.

A \textit{trivial loop\/} in a graph is a length one cycle which bounds a disk face.
By $(*)$, neither $G_S$ nor $G_T$ contains trivial loops.

A family of edges $\{e_1,e_2,\dots,e_p\}$ in $G_S$
is a \textit{Scharlemann cycle\/} (of length $p$)
if it bounds a disk face of $G_S$, and all the edges have the same pair
of labels $\{x,x+1\}$ at their two endpoints, which is called
the \textit{label pair\/} of the Scharlemann cycle.
Note that each edge $e_i$ connects the vertex $u_x$ with $u_{x+1}$ in $G_T$.
A Scharlemann cycle of length two is called an \textit{$S$-cycle} for short.
Remark that the interior of the face bounded by a Scharlemann cycle may
meet $\widehat{T}$, since $T$ is not necessarily incompressible in $E(K)$.

Let $\sigma$ be a Scharlemann cycle in $G_S$ with label pair $\{x,x+1\}$.
If the edges of $\sigma$ (and vertices $u_x$ and $u_{x+1}$) are contained
in an essential annulus $A$ in $\widehat{T}$, 
and if they do not lie in a disk in $\widehat{T}$, then
we say that the edges of $\sigma$ \textit{lie in an essential annulus} in $\widehat{T}$.

\begin{lemma}  
Let $\sigma$ be a Scharlemann cycle in $G_S$ of length $p$
with label pair $\{x,x+1\}$, where $p$ is $2$ or $3$.
Let $f$ be the face of $G_S$ bounded by $\sigma$.
If the edges of $\sigma$ do not lie in a disk in $\widehat{T}$, then
they lie in an essential annulus $A$ in $\widehat{T}$.
Furthermore, if $\mathrm{Int}f\cap \widehat{T}=\emptyset$, then
$M=N(A\cup H_{x,x+1}\cup f)$ is a solid torus such that
the core of $A$ runs $p$ times in the longitudinal direction of $M$. 
\end{lemma}

\begin{proof}
If $p=2$, then it is obvious that
the edges of $\sigma$ lie in an essential annulus in $\widehat{T}$.

Assume $p=3$.
Let $\sigma=\{e_1,e_2,e_3\}$.
If the endpoints of $e_1,e_2,e_3$ appear in this order when one travels
around $u_x$ clockwise, say, then those of $e_1,e_2,e_3$ appear in the same order
when one travels around $u_{x+1}$ anticlockwise, since $u_x$ and $u_{x+1}$ have
distinct parities.
This observation implies that the edges of $\sigma$ lie
in an essential annulus in $\widehat{T}$.

Consider the genus two handlebody $N(A\cup H_{x,x+1})$.
Then $M$ is obtained by attaching a $2$-handle $N(f)$.
Since there is a meridian disk of $N(A)$ which intersects $\partial f$ once,
$\partial f$ is primitive and therefore $M$ is a solid torus.
It is not hard to see that the core of $A$ runs $p$ times in the longitudinal
direction of $M$.  See also \cite[Lemma 3.7]{GL2}
\end{proof}

\begin{lemma}  
Let $\xi$ be a loop in $S\cap T$.
Suppose that $\xi$ bounds a disk $\delta$ in $S$ with
$\mathrm{Int}\delta\cap\widehat{T}=\emptyset$.
If $\xi$ is inessential in $\widehat{T}$,
then all vertices of $G_T$ must lie in the disk bounded by $\xi$.
\end{lemma}

\begin{proof}
Let $\delta'$ be the disk bounded by $\xi$ in $\widehat{T}$.
Then $\delta'\cap V\ne\emptyset$, since $\xi$ is essential in $T$ by the assumption
on $S\cap T$.
If both sides of $\xi$ on $\widehat{T}$ meet $V$,
replace $\widehat{T}$ by $\widehat{T'}=(\widehat{T}-\delta')\cup\delta$.
Then $\widehat{T'}$ gives a new Heegaad torus of $K(r)$ satisfying $(*)$.
However this contradicts the choice of $\widehat{T}$,
since $|\widehat{T'}\cap K^*|<|\widehat{T}\cap K^*|$.
Hence all vertices of $G_T$ lie in $\delta'$.
\end{proof}

\begin{lemma}  
Let $\sigma$ be a Scharlemann cycle in $G_S$ of length $p$ with
label pair $\{x,x+1\}$,
and let $f$ be the face of $G_S$ bounded by $\sigma$.
Suppose that $p\ne r$.
Then the edges of $\sigma$ cannot lie in a disk in $\widehat{T}$, and
$\mathrm{Int}f\cap \widehat{T}=\emptyset$.
\end{lemma}

\begin{proof}
Assume for contradiction that the edges of $\sigma$ lie in a disk $D$ in $\widehat{T}$.
Let $\Gamma$ be the subgraph of $G_T$ consisting of two vertices $u_x$ and $u_{x+1}$
along with the edges of $\sigma$.

First, suppose that $\mathrm{Int}f\cap D\ne\emptyset$.
By the cut-and-paste operation of $f$, it can be assumed that
any component in $\mathrm{Int}f\cap D$ is essential in $D-\Gamma$.
Therefore all components in $\mathrm{Int}f\cap D$ are parallel to $\partial D$ in $D-\Gamma$. 
Then we can replace $D$ by a subdisk which does not meet
$\mathrm{Int}f$.
We may now assume that $\mathrm{Int}f\cap D=\emptyset$.
Then $N(D\cup H_{x,x+1}\cup f)$ gives a punctured lens space.
Since a lens space $K(r)$ is irreducible, this means that $K(r)$ is a lens space
whose fundamental group has order $p$.
This contradicts the assumption that $p\ne r$.
Thus the edges of $\sigma$ cannot lie in a disk in $\widehat{T}$.

Assume that $\mathrm{Int}f\cap \widehat{T}\ne\emptyset$.
Let $\mu$ be an innermost component of $\mathrm{Int}f\cap\widehat{T}$ on $f$. 
By Lemma 2.2, $\mu$ is essential in $\widehat{T}$.
Then it can be assumed that the disk $\delta$ bounded by $\mu$ on $f$ is contained in $W$,
say, one of the solid tori bounded by $\widehat{T}$ in $K(r)$.
Thus $\delta$ is a meridian disk of $W$.

In $W$, compress $\widehat{T}$ along $\delta$ to obtain a $2$-sphere $Q$.
There is a disk $E$ in $Q$ which contains the edges of $\sigma$ and two vertices
$u_x$ and $u_{x+1}$.
Even if $\mathrm{Int}f\cap E\ne\emptyset$, the cut-and-paste operation gives
a new $f$ with $\mathrm{Int}f\cap E=\emptyset$.
Thus $N(E\cup H_{x,x+1}\cup f)$ gives a punctured lens space whose
fundamental group has order $p$, which contradicts the assumption again.
\end{proof}

When there exist two Scharlemann cycles with disjoint label pairs,
the assumption on the length in the statement of Lemma 2.3 is not necessary.

\begin{lemma} 
Let $\sigma_1$ and $\sigma_2$ be Scharlemann cycles in $G_S$ with disjoint label pairs,
and let $f_1$ and $f_2$ be the faces of $G_S$ bounded by $\sigma_1$ and $\sigma_2$ respectively.
Then the edges of $\sigma_i$ lie in an essential annulus $A_i$ in $\widehat{T}$ with
$A_1\cap A_2=\emptyset$,
and $\mathrm{Int}f_i\cap\widehat{T}=\emptyset$ for $i=1,2$.
\end{lemma}

\begin{proof}
Let $\{x_i,x_i+1\}$ be the label pair of $\sigma_i$.
Assume that the edges of $\sigma_1$ lie in a disk $D_1$ in $\widehat{T}$
for contradiction.
By the same argument in the proof of Lemma 2.3,
we may assume that $\mathrm{Int}f_1\cap D_1=\emptyset$.
If $\mathrm{Int}f_1\cap\widehat{T}=\emptyset$, then 
$N(D_1\cup H_{x_1,x_1+1}\cup f_1)$ gives a punctured lens space in a solid torus,
which is impossible.  Therefore $\mathrm{Int}f_1\cap\widehat{T}\ne\emptyset$.

Choose an innermost component $\xi$ of
$\mathrm{Int}f_1\cap\widehat{T}$ on $f_1$.
Let $\delta$ be the disk bounded by $\xi$ on $f_1$.

Assume that $\xi$ is inessential in $\widehat{T}$. By Lemma 2.2, 
$G_T$ lies in the disk bounded by $\xi$.
Then the edges of $\sigma_2$ also lie in a disk $D_2$ in $\widehat{T}$.
We remark that one of $D_1$ and $D_2$ may be contained in the other, possibly. 
As above, we can assume that $\mathrm{Int}f_2\cap D_2=\emptyset$.

If $D_1\cap D_2=\emptyset$, then
we can assume that $\mathrm{Int}f_i\cap D_j=\emptyset$ for $i,j\in\{1,2\}$
by the cut-and-paste operation of $f_i$.

Otherwise, $D_2\subset D_1$, say.
Clearly, $\mathrm{Int}f_1\cap D_j=\emptyset$ for $j=1,2$.
If $\mathrm{Int}f_2\cap D_1\ne\emptyset$, then it can be assumed that each component of
$\mathrm{Int}f_2\cap D_1$ is parallel to $\partial D_2$ in $D_1-\Gamma$,
where $\Gamma$ is the subgraph of $G_T$, consisting of
the vertices $u_{x_i}$,$u_{x_i+1}$ along with the edges of $\sigma_i$ for $i=1,2$.
But this contradicts Lemma 2.2.
Therefore, we can assume that $\mathrm{Int}f_i\cap D_j=\emptyset$ for $i,j\in\{1,2\}$
in either case.

Then $N(D_1\cup H_{x_1,x_1+1}\cup f_1)$ and $N(D_2\cup H_{x_2,x_2+1}\cup f_2)$
give two disjoint punctured lens spaces in $K(r)$, which is impossible.
(When $D_2\subset D_1$, say, we have to push $D_2$ into a suitable direction
away from $D_1$.)

Therefore $\xi$ is essential in $\widehat{T}$.
Then $\delta$ is a meridian disk of the solid torus $W$, say.
Compressing $\widehat{T}$ along $\delta$ gives a $2$-sphere $Q$ on which
there are two disjoint disks $E_1,E_2$ each containing the edges of $\sigma_1,\sigma_2$,
respectively.  Then the same argument as above gives a contradiction.

Therefore the edges of $\sigma_i$ cannot lie in a disk in $\widehat{T}$ for $i=1,2$,
and then
there are disjoint essential annuli $A_i$ in $\widehat{T}$ in which
the edges of $\sigma_i$ lie for $i=1,2$, respectively.

Suppose that $\mathrm{Int}f_1\cap\widehat{T}\ne\emptyset$.
Consider an innermost component $\eta$ of $\mathrm{Int}f_1\cap\widehat{T}$
in $f_1$. 
By Lemma 2.2, $\eta$ is essential in $\widehat{T}$.
As above, there are two disjoint punctured lens spaces in $K(r)$, which is impossible again.
Similarly for $f_2$.
Therefore $\mathrm{Int}f_i\cap\widehat{T}=\emptyset$ for $i=1,2$.
\end{proof}

Let $f$ be a face of $G_S$.
Although $\mathrm{Int}f\cap \widehat{T}\ne\emptyset$ in general,
a small collar neighborhood of $\partial f$ in $f$ is contained in one side of $\widehat{T}$.
Then we say that $f$ \textit{lies on that side of} $\widehat{T}$.

The next two lemmas deal with the situation where $G_S$ has two Scharlemann cycles
of length two and three simultaneously.

\begin{lemma}  
Let $\sigma$ be an $S$-cycle in $G_S$,
and let $\tau$ be a Scharlemann cycle in $G_S$ of length three.
Let $f$ and $g$ be the faces of $G_S$ bounded by $\sigma$ and $\tau$ respectively.
If $\sigma$ and $\tau$ have disjoint label pairs, then
$\sigma$ and $\tau$ lie on opposite sides of $\widehat{T}$, and
$r\not\equiv 0 \pmod{3}$.
\end{lemma}

\begin{proof}
Let $\{x,x+1\}$, $\{y,y+1\}$ be the label pairs of $\sigma$ and $\tau$, respectively.
By Lemma 2.4, the edges of $\sigma$ give an essential cycle in $\widehat{T}$
after shrinking two fat vertices $u_x$ and $u_{x+1}$ to points, and
$\mathrm{Int}f\cap \widehat{T}=\emptyset$.
Then $f$ is contained in the solid torus, $W$ say, and
the union $H_{x,x+1}\cup f$ gives a M\"obius band $B$ properly embedded in $W$, after
shrinking $H_{x,x+1}$ to its core radially.  See Figure 1.


\bigskip
$$\fbox{Figure 1}$$
\bigskip

Similarly, by Lemma 2.4, the edges of $\tau$
lie in an essential annulus $A$ in $\widehat{T}$ which is disjoint from
the edges of $\sigma$, and $\mathrm{Int}g\cap \widehat{T}=\emptyset$.

Suppose that $g\subset W$. 
If a solid torus $J$ is attached to $W$ along their boundaries so that
the slope of $\partial B$ bounds a meridian disk of $J$, then
the resulting manifold $N=J\cup W$ contains a projective plane, and therefore
$N=L(2,1)$.
However, in $N$, the edges of $\tau$ are contained in a disk $D$ obtained by
capping a boundary component of $A$ off by a meridian disk of $J$.
Then $N(D\cup H_{y,y+1}\cup g)$ gives a punctured lens space of order three in $N$, which is
impossible.
Thus $f$ and $g$ lie on opposite sides of $\widehat{T}$.

Next, assume that $r\equiv 0 \pmod{3}$ for contradiction.
We may assume that $f\subset W$ and $g\subset U$.

By Lemma 2.1, $M=N(A\cup H_{y,y+1}\cup g)$ is a solid 
torus, and $A$ runs three times in the longitudinal 
direction on $\partial M$. 
The annulus $A'=\mathrm{cl}(\partial M-A)$ is 
properly embedded in $U$, and so $A'$ is parallel to 
$\mathrm{cl}(\widehat{T}-A)$. Therefore, $A$ runs three times 
in the longitudinal direction of $U$.

The slope determined by $\partial M$ on $\partial W$ meets
a meridian of $W$ twice.
On $\partial U$, the slope can be expressed $a/3$
and the meridian of $W$ defines a slope $b/r$ for some integers $a, b$.
Then $\Delta(a/3,b/r)=|ar-3b|=3|ar/3-b|\ne 2$, which is a contradiction.
\end{proof}

\begin{lemma}  
Let $\sigma$, $\tau$, $f$, $g$ be as in Lemma 2.5.
Suppose that $\sigma$ and $\tau$ lie on opposite sides of $\widehat{T}$ and 
have the same label pair, and that $r\ne 2,3$.
If there is an essential annulus $A$ in $\widehat{T}$ in which
the edges of $\sigma$ and $\tau$ lie, then $r\not\equiv 0 \pmod{3}$.
\end{lemma}

\begin{proof}
By Lemma 2.3, $\mathrm{Int}f\cap \widehat{T}=\emptyset$
and $\mathrm{Int}g\cap \widehat{T}=\emptyset$.
We remark that $t=2$.
Hence $\sigma$ and $\tau$ have the label pair $\{1,2\}$.

We may assume that
$H_{1,2}\subset W$ and $f\subset W$.
Then $M_1=N(A\cup H_{1,2}\cup f)$ is a solid torus, and $A$ runs twice in the longitudinal
direction on $\partial M_1$ by Lemma 2.1.
Furthermore, the annulus $A'_1=\mathrm{cl}(\partial M_1-A)$ is parallel
to $\mathrm{cl}(\widehat{T}-A)$ in $W$.
Similarly, $M_2=N(A\cup H_{2,1}\cup g)$ is a solid torus,
and $A$ runs three times in the longitudinal direction on $\partial M_2$ by Lemma 2.1.
The annulus $A'_2=\mathrm{cl}(\partial M_2-A)$ is also parallel
to $\mathrm{cl}(\widehat{T}-A)$ in $U$.
Then the same argument as in the proof of Lemma 2.5 gives the desired result.
\end{proof}

\section{The generic case}
In this section we prove Theorem 1.1 under the hypothesis $t\ge 4$.
The case $t=2$ will be dealt with separately in the next section.

\begin{lemma}  
If $G_S$ has a family of more than $t$ mutually parallel edges, then
there are at least two $S$-cycles on disjoint label pairs in the family.
\end{lemma}

\begin{proof}
Assume that there are $t+1$ mutually parallel edges $e_1,e_2,\dots,e_{t+1}$ in
$G_S$, numbered  successively.
We may assume that $e_i$ has the label $i$ at one endpoint for $1\le i\leq t$,
and $e_{t+1}$ has the label $1$.
By the parity rule, $e_{2j}$ has the label $1$ at the other endpoint for some $j$.
If $2j=t$, then $\{e_j,e_{j+1}\}$ and $\{e_t,e_{t+1}\}$ form $S$-cycles
with disjoint label pairs.
If $2j<t$, then $\{e_j,e_{j+1}\}$ and
$\{e_{t/2+j},e_{t/2+j+1}\}$ form $S$-cycles with disjoint label pairs.
\end{proof}

\begin{lemma}  
If $G_S$ contains two $S$-cycles on disjoint label pairs, then $K(r)$ is
$L(4k,2k\pm 1)$ for some $k\ge 1$.
\end{lemma}

\begin{proof}
Let $\sigma_1$ and $\sigma_2$ be $S$-cycles in $G_S$ with the label pairs
$\{x_1,x_1+1\},\{x_2,x_2+1\}$, respectively, where $\{x_1,x_1+1\}\cap\{x_2,x_2+1\}=\emptyset$.
Let $f_i$ be the face of $G_S$ bounded by $\sigma_i$.
By Lemma 2.4, the edges of $\sigma_i$ lie in an essential annulus in $\widehat{T}$, and
$\mathrm{Int}f_i\cap \widehat{T}=\emptyset$.
If we shrink $H_{x_i,x_i+1}$ to its core radially, then $H_{x_i,x_i+1}\cup \sigma_i$
gives a M\"obius band $B_i$ properly embedded in $U$ or $W$.

Since $\partial B_i$ is essential in $\widehat{T}$,
$\partial B_1$ and $\partial B_2$ are parallel in $\widehat{T}$.
Thus the union of $B_1\cup B_2$ and an annulus in $\widehat{T}$
bounded by $\partial B_1$ and $\partial B_2$ gives a Klein bottle in $K(r)$.
(In fact, $f_1$ and $f_2$ lie on opposite sides of $\widehat{T}$, otherwise
a Klein bottle will be found in a solid torus.)
It is known that a lens space contains a Klein bottle if and only if
the lens space has the form of
$L(4k,2k\pm 1)$ for some $k\ge 1$ \cite[Corollary 6.4]{BW}.
\end{proof}

It is conjectured that $L(4,1)$ cannot be obtained from non-trivial knots in $S^3$
by Dehn surgery \cite{BL,Go01}.
In general, it seems to be unknown that $L(4k,2k\pm 1)$ can arise by surgery
on hyperbolic knots (we expect that it cannot).

\begin{lemma} 
Let $\{e_1,e_2,\dots,e_t\}$ be mutually parallel edges in $G_S$ numbered successively.
Then $\{e_{t/2},e_{t/2+1}\}$ is an $S$-cycle.
\end{lemma}

\begin{proof}
We may assume that $e_i$ has the label $i$ at one endpoint for $1\le i\le t$.
If $e_t$ has the label $1$ at the other endpoint, then
$\{e_{t/2},e_{t/2+1}\}$ is an $S$-cycle.
Therefore we suppose that $e_{2j}$ has the label $1$ at the other endpoint for some $j<t/2$.
Then $\sigma_1=\{e_j,e_{j+1}\}$ and
$\sigma_2=\{e_{t/2+j},e_{t/2+j+1}\}$
form $S$-cycles with disjoint label pairs.
By the same argument as in the proof of Lemma 3.2,
we obtain two disjoint M\"obius bands $B_1, B_2$ from the faces of $\sigma_1,\sigma_2$,
and furthermore a Klein bottle $F$ in $K(r)$.
(If $t/2$ is even, then these M\"obius bands lie on the same side of $\widehat{T}$ and
therefore $F$ is contained in a solid torus, which is impossible.
Hence we have $t/2$ is odd.)
Since $\partial B_1$ and $\partial B_2$ are parallel in $\widehat{T}$,
they divide $\widehat{T}$ into two annuli $A_1$ and $A_2$.

In $G_T$, $u_k$ and $u_{2j-k+1}$ lie in the
same annulus for $1\le k\le j-1$, since the edge $e_k$ connects the two vertices in $G_T$.
Similarly,
$u_{2j+\ell}$ and $u_{t+1-\ell}$ for $1\le \ell\le t/2-j-1$
lie in the same annulus.
Therefore we see that $\mathrm{Int}A_i$ contains an even number of vertices for $i=1,2$.
We may assume that $F$ is obtained as the union $B_1\cup B_2\cup A_1$.
However, $F$ meets $K^*$ in an even number of points (after a perturbation).
Then $F'=F\cap E(K)$ gives a punctured Klein bottle properly embedded in $E(K)$
having an even number of boundary components.
By attaching suitable annuli in $\partial E(K)$ to $F'$ along boundaries,
we have a closed non-orientable surface in $E(K)$, which is impossible.
\end{proof}

The reduced graph $\overline{G}_S$ of $G_S$ is defined to be the graph obtained
from $G_S$ by amalgamating each set of mutually parallel edges of $G_S$ to
a single edge.
If an edge $\overline{e}$ of $\overline{G}_S$ corresponds to $s$ mutually parallel edges
of $G_S$, then the \textit{weight\/} of $\overline{e}$ is defined to be $s$, and
we denote by $w(\overline{e})=s$.
If $w(\overline{e})=t$, then $e$ is called a \textit{full edge\/}.

\begin{proposition}  
If $t\ge 4$, then $r\le 12g-7$, where $g$ is the genus of $K$.
\end{proposition}

\begin{proof}
Since $G_S$ does not contain trivial loops, the unique vertex $v$ has 
valency at most $12g-6$ in $\overline{G}_S$ (see \cite[Lemma 6.2]{GLi}).
Therefore the edges of $G_S$ are partitioned into at most $6g-3$ families of
parallel edges.

If there is an edge $\overline{e}$ in $\overline{G}_S$ with $w(\overline{e})>t$, then
$K(r)=L(4k,2k\pm 1)$ for some $k\ge 1$ by Lemmas 3.1 and 3.2, that is,
$K(r)$ contains a Klein bottle.
Then $r\le 12g-8$ \cite{Te3}.

Hence suppose that $w(\overline{e})\le t$ for any edge $\overline{e}$ of $\overline{G}_S$.
Recall that the vertex $v$ has valency $rt$ in $G_S$.
Then $rt\le (12g-6)t$, hence $r\le 12g-6$.

Finally, suppose that $r=12g-6$.
Then any edge of $\overline{G}_S$ is full, and each face of $\overline{G}_S$ is a $3$-sided disk.
By Lemma 3.3,
we may assume that $G_S$ contains an $S$-cycle with label pair $\{t/2,t/2+1\}$
and a Scharlemann cycle of length three with label pair $\{t,1\}$.
Then $r\not\equiv 0 \pmod{3}$ by Lemma 2.5, which is a contradiction.
Therefore $r\le 12g-7$.
\end{proof}

\section{The case that $t=2$}

By the parity rule, each edge of $G_T$ connects different vertices $u_1$ and $u_2$.
Then there are four \textit{edge classes\/} in $G_T$, i.e., isotopy classes
of non-loop edges of $G_T$ in $\widehat{T}$ rel $u_1\cup u_2$.
They are called $1,\alpha,\beta,\alpha\beta$ as illustrated in Figure 2
(see \cite[Figure 7.1]{GL2}).


\bigskip
$$\fbox{Figure 2}$$
\bigskip

We label an edge of $e$ of $G_S$ by the class of the corresponding edge of $G_T$,
and we call the label the \textit{edge class label\/} of $e$.

For a face $f$ of $G_S$, if a small collar neighborhood of $\partial f$ in $f$
is contained in $U$\ ($W$) ,
then $f$ is said to be \textit{black\/} (resp. \textit{white}).

\begin{lemma}  
Suppose that $r\ne 2$. 
Then any two black (white) bigons in $G_S$ have the same pair of edge class labels.
\end{lemma}

\begin{proof}
By Lemma 2.3, the interior of a black (white) bigon is disjoint from $\widehat{T}$. 
Then the proof of \cite[Lemma 5.2]{GL3} remains valid.
Remark that a final contradiction comes from the fact that
a Klein bottle will be found in a solid torus $U$ or $W$. 
\end{proof}

\begin{lemma}   
Let $e$ and $e'$ be edges of $G_S$.
If $e$ and $e'$ are parallel in $G_S$, then they have distinct edge class labels.
\end{lemma}

\begin{proof}
If $e$ and $e'$ are parallel in $G_S$ and have the same edge class label,
then they are also parallel in $G_T$. Then $E(K)$ contains
a M\"obius band by \cite[Lemma 2.1]{Go2}, which contradicts the fact that $K$ is hyperbolic.
\end{proof}

\begin{lemma}  
If $r\ne 2$, then
$G_S$ cannot contain more than 3 mutually parallel edges.
\end{lemma}

\begin{proof}
Suppose that there are $4$ mutually parallel edges.
Then there are two bigons with the same color among these $4$ parallel edges.
By Lemma 4.1, these two bigons have the same pair of edge class labels.
This contradicts Lemma 4.2. 
\end{proof}

\begin{lemma}       
Suppose that $r\ne 2$. 
If $G_S$ contains a black bigon and a white bigon which have an edge
in common, then the other faces of $G_S$ are not bigons.
\end{lemma}

\begin{proof}
Let $e_1,e_2,e_3$ be adjacent parallel edges of $G_S$.
By Lemma 4.2, these three edges have distinct edge class labels.
Let $\lambda,\mu,\nu$ be the edge class labels of $e_1,e_2,e_3$ respectively.
Let us denote the endpoints of $e_i$ by $\partial^j e_i$ for $j=1,2$.
See Figure 3.


\bigskip
$$\fbox{Figure 3}$$
\bigskip

Note that $\partial^1e_1$ and $\partial^1e_3$ appear consecutively around
the vertex $u_1$ in the order, when traveling around $\partial u_1$ anticlockwise, say.  
Then $\partial^2e_3$ and $\partial^2e_1$ appear consecutively around $u_2$ in the order,
when traveling around $\partial u_2$ clockwise.
These come from the facts that $r$ is integral, and that $u_1$ and $u_2$ have
distinct parities.
Then there is no other edge of edge class $\lambda$\ $(\nu)$ than $e_1$\ $(e_3)$ in $G_T$.
The conclusion follows from Lemma 4.1.
\end{proof}

\begin{proposition} 
If $t=2$, then $r\le 12g-7$.
\end{proposition}

\begin{proof}
The unique vertex $v$ has 
valency at most $12g-6$ in $\overline{G}_S$, and
the edges of $G_S$ are partitioned into at most $6g-3$ families of
parallel edges.
Recall that $v$ has valency $2r$ in $G_S$.

By Lemma 4.3, $G_S$ cannot contain $4$ mutually parallel edges.
If $G_S$ contains $3$ mutually parallel edges, then
we have $r\le (6g-3)+2=6g-1$ by Lemma 4.4.

If $G_S$ does not contain $3$ mutually parallel edges,
then each edge of $\overline{G}_S$ has weight $1$ or $2$.
Hence $r\le 2(6g-3)=12g-6$.

Suppose that $r=12g-6$.
Then any edge of $\overline{G}_S$ is full, and hence
$G_S$ has $6g-3$ black, say, bigons and each white face of $G_S$ is a $3$-sided disk.
Therefore there are an $S$-cycle $\sigma$ and a Scharlemann cycle $\tau$ of
length three in $G_S$ with the same label pair $\{1,2\}$.
By Lemma 4.1, all black bigons have the same pair of edge class label $\{\lambda,\mu\}$, say.
Then the edges of $\tau$ have the same edge class labels $\lambda,\mu$ by Lemma 2.3.
This means that there is an essential annulus $A$ in $\widehat{T}$ which contains
the edges of $\sigma$ and $\tau$.
By Lemma 2.6, we have $r\not\equiv 0\pmod{3}$, which is a contradiction.
Therefore $r\le 12g-7$.
\end{proof}

\bigskip
\noindent
\textit{Proof of Theorem 1.1.}
This follows immediately from Propositions 3.4 and 4.5.
\hfill\qedsymbol

\section{Genus one case: the case $t\ge 4$}

In the remainder of this paper, $K$ is assumed to be a genus one, hyperbolic knot in $S^3$
in order to prove Theorem 1.2.
First, we deal with the case $t\ge 4$ in this section.

\begin{theorem} 
If $K$ has genus one, then $K(r)$ is not $L(2,1), L(4k,2k\pm 1)$ for any $k\ge 1$.
\end{theorem}

\begin{proof}
This follows from \cite{Do,Te1,Te2}.
\end{proof}

\begin{lemma} 
$G_S$ cannot contain two $S$-cycles on disjoint label pairs.
\end{lemma}

\begin{proof}
If $G_S$ contains two $S$-cycles on disjoint label pairs,
then $K(r)$ is $L(4k,2k\pm 1)$ for some $k\ge 1$ by Lemma 3.2.
But this is impossible by Theorem 5.1.
\end{proof}

\begin{lemma}  
If $r$ is odd, then
$G_S$ cannot have more than $t/2$ mutually parallel edges.
\end{lemma}

\begin{proof}
The vertex $v$ has valency $rt$ in $G_S$.
Recall that the edges of $G_S$ are partitioned into at most three
families of mutually parallel edges.
Let $A$ be a family of mutually parallel edges in $G_S$, and
suppose that $A$ consists of more than $t/2$ edges, 
$a_1,a_2,\dots,a_p$ numbered consecutively.
Note that $p\le t$ by Lemmas 3.1 and 5.2.
We may assume that $a_i$ has the label $i$ at one endpoint for $1\le i\le p$.
Then $a_p$ has the label $t/2+1$ at the other endpoint, since $r$ is odd.
See Figure 4.


\bigskip
$$\fbox{Figure 4}$$
\bigskip

By the parity rule, $p\ne t/2+1$.
Thus $p>t/2+1$.
Then $\{a_{(t/2+p)/2},a_{(t/2+p)/2+1}\}$ forms an $S$-cycle.
Furthermore, some edge between $a_2$ and $a_{t/2}$ has the label $1$ at the other endpoint.
Therefore, there is another $S$-cycle whose label pair is disjoint from
that of the above $S$-cycle.
This contradicts Lemma 5.2.
\end{proof}

By Proposition 3.4, we have that $r\le 5$.
In fact, the cases that $r=3,5$ remain by Theorem 5.1.

\begin{lemma} 
The case that $r=3$ is impossible.
\end{lemma}

\begin{proof}
The vertex $v$ has valency $3t$ in $G_S$.
By Lemma 5.3, $G_S$ consists of three families of mutually parallel edges, each
containing exactly $t/2$ edges.
Then there is no $S$-cycle in $G_S$, but there are two Scharlemann cycles
$\tau_1$ and $\tau_2$ of length three in $G_S$.
Let $g_i$ be the face of $G_S$ bounded by $\tau_i$ for $i=1,2$.
We may assume that $g_1$ has the label pair $\{t,1\}$, and $g_2$ has $\{t/2,t/2+1\}$.

By Lemma 2.4, 
there are disjoint essential annuli $A_i$ in $\widehat{T}$ in which
the edges of $\tau_i$ lie, and $\mathrm{Int}g_i\cap\widehat{T}=\emptyset$ for $i=1,2$.

\begin{claim}  
The faces $g_1$ and $g_2$ lie on opposite sides of $\widehat{T}$.
\end{claim}

\noindent
\textit{Proof of Claim 5.5.}
Suppose that $g_i\subset W$, say,  for $i=1,2$.
Let $s$ be the slope on $\partial W$ determined by the essential annuli $A_i$.
Performing $s$-Dehn filling on $W$, that is, attaching a solid torus $J$
to $W$ along their boundaries so that $s$ bounds a meridian disk of $J$,
we obtain a closed $3$-manifold $M$, which is either $S^3$, $S^2\times S^1$ or a lens space.
However, there are
two disjoint disks $D_1$ and $D_2$, which contain the edges of $\tau_1$ and $\tau_2$,
respectively, on the $2$-sphere $Q$ obtained by compressing $\widehat{T}$
along $s$ by a meridian disk of $J$.
Then $N(D_1\cup H_{t,1}\cup g_1)$ and $N(D_2\cup H_{t/2,t/2+1}\cup g_2)$ give
two punctured lens spaces in $M$,
which is impossible.
\hfill (Proof of Claim 5.5)\qedsymbol
\bigskip

Therefore, $t/2$ and $t$ must have opposite parities, and
so $t/2$ is odd.
In particular, $t\ge 6$.

In $G_S$, there are exactly three edges whose endpoints have the pair of labels
$\{j,t+1-j\}$ for $j=1,2,\dots,t/2$.
Therefore, $G_T$ consists of $t/2$ components, each consisting two vertices
$u_j$ and $u_{t+1-j}$ along with three edges connecting them.

\begin{claim} 
Each component of $G_T$ does not lie in a disk in $\widehat{T}$.
\end{claim}

\noindent
\textit{Proof of Claim 5.6.}
If there is a component of $G_T$ which lie in a disk in $\widehat{T}$,
then we can take an innermost one $\Lambda$.
That is, $\Lambda$ lies in a disk $D$ in $\widehat{T}$, and
there is no other component of $G_T$ in $D$.

Consider the intersection between $D$ and $S$.
Then $S$ is divided into two disks $g_3$ and $g_4$ by the edges of $\Lambda$.
By the cut-and-paste operation of $g_3$ or $g_4$, and taking $D$ by a smaller
one, if necessary, we can assume that $\mathrm{Int}g_3$ and $\mathrm{Int}g_4$
do not meet $D$. 
Then $N(D\cup V\cup g_3\cup g_4)$, where $V$ is the attached solid torus, gives
a connected sum of two lens spaces minus an open $3$-ball in $K(r)$, which
is impossible.
\begin{flushright}
(Proof of Claim 5.6)\qedsymbol
\end{flushright}
\bigskip

Thus, we may assume that $A_i$ contains only the edges and vertices of $\tau_i$ for $i=1,2$.

Assume that $g_1\subset W$ and $g_2\subset U$.
Let $M_1=N(A_1\cup H_{t,1}\cup g_1)$ and $M_2=N(A_2\cup H_{t/2,t/2+1}\cup g_2)$.
Let $A'_i=\mathrm{cl}(\partial M_i-A_i)$ for $i=1,2$.
Then $A'_1$ is a properly embedded annulus in $W$ and $A'_2$ is a
properly embedded annulus in $U$.
By Lemma 2.1, $M_i$ is a solid torus such that the core of $A_i$ runs
three times in the longitudinal direction of $M_i$ for $i=1,2$.
Therefore, $A'_1$ is parallel to the annulus $\mathrm{cl}(\widehat{T}-A_1)$ in $W$,
and $A'_2$ is parallel to $\mathrm{cl}(\widehat{T}-A_2)$ in $U$.
Let $\widehat{T'}=(\widehat{T}-(A_1\cup A_2))\cup A'_1\cup A'_2$.
Then it is easy to see that $\widehat{T'}$ is a new Heegaard torus in $K(r)$
such that $|\widehat{T'}\cap V|=t-4 \ (>0)$.
Furthermore, $\widehat{T'}$ satisfies $(*)$, which contradicts the choice of $\widehat{T}$.
\end{proof}


\begin{lemma}       
The case that $r=5$ is impossible.
\end{lemma}

\begin{proof}
Since the vertex $v$ has valency $5t$ in $G_S$, there are
more than $t/2$ mutually parallel edges in $G_S$, which contradicts Lemma 5.3.
\end{proof}

\bigskip
\noindent
\textit{Proof of Theorem 1.2 when $t\ge 4$.}
By Proposition 3.4, $r\le 5$, and in fact, the remaining cases are $r=3,5$ by Theorem 5.1.
But these cases are impossible by Lemmas 5.4 and 5.7.
\hfill\qedsymbol

\section{Genus one case: the case $t=2$}

In the case that $t=2$, the following lemma plays a key role.

Recall that an \textit{unknotting tunnel\/} $\gamma$ for a knot or link $K$ in $S^3$ is
a simple arc properly embedded in the exterior $E(K)$ such that
$\mathrm{cl}(E(K)-N(\gamma))$ is homeomorphic to a handlebody of genus two.

\begin{lemma} 
Let $K$ be a genus one knot in $S^3$, and let $S$ be a minimal
genus Seifert surface of $K$.
If $K$ has an unknotting tunnel $\gamma$ such that $\gamma\subset S$, then
$K$ is $2$-bridge.
\end{lemma}

\begin{proof}
Take a regular neighborhood $N$ of $\gamma$ in $S$.
Let $F=\mathrm{cl}(S-N)$.
Then $F$ is an annulus whose boundary defines a link $L$ in $S^3$.
Note that $F$ is incompressible in the exterior of $L$, and
$L$ has an unknotting tunnel.
Then $L$ is a $2$-bridge torus link by \cite[Theorem 1]{EU}.
Furthermore, an unknotting tunnel of such a link is determined by
\cite{AR}.
Then $S$ can be restored  by taking the union of $F$ and $N$, showing that
$K$ is $2$-bridge.
\end{proof}

\begin{lemma} 
Suppose that $K$ has genus one.
Then $G_S$ cannot have more than two mutually parallel edges.
\end{lemma}

\begin{proof}
If there are three mutually parallel edges in $G_S$,
there are two $S$-cycles $\sigma_1$ and $\sigma_2$ whose faces $f_1$ and $f_2$
lie on opposite sides of $\widehat{T}$.
Since $K(r)\ne L(2,1)$ by Theorem 5.1,
we can assume that $\mathrm{Int}f_i\cap \widehat{T}=\emptyset$ for $i=1,2$ by Lemma 2.3.
Then we may also assume that $f_1, H_{1,2}\subset W$ and $f_2, H_{2,1}\subset U$.
Note that $\mathrm{cl}(W-H_{1,2})$ and $\mathrm{cl}(U-H_{2,1})$ are
handlebodies of genus two, since the core of $H_{1,2}$ ($H_{2,1}$) lies
on a M\"obius band which is obtained from $H_{1,2}\cup f_{1}$
($H_{2,1}\cup f_{2}$) by shrinking $H_{1,2}$ (resp. $H_{2,1}$) to its core radially.

Let $\alpha$ and $\beta$ be the arc components of $f_1\cap H_{1,2}$.
Let $\gamma$ be a simple arc in $f_1$ which connects a point in $\alpha$
with one in $\beta$.  See Figure 5.


\bigskip
$$\fbox{Figure 5}$$
\bigskip

Then it can be seen that $\mathrm{cl}(W-H_{1,2}-N(\gamma))$ is homeomorphic to
$T\times I$, where $I$ denotes an interval.
Therefore, $\gamma$ gives an unknotting tunnel of $K$ which lies on $S$. 
By Lemma 6.1, $K$ is $2$-bridge, which contradicts the fact that
a hyperbolic $2$-bridge knot has no cyclic surgery \cite{Ta}.
\end{proof}

The remaining cases are $r=3,5$ again by Proposition 4.5 and Theorem 5.1.

\begin{lemma} 
The case $r=3$ is impossible.
\end{lemma}

\begin{proof}
Recall that $\widehat{T}$ is separating in $K(r)$, and therefore the faces of $G_S$ are
partitioned into black and and white ones.
This implies that $G_S$ has no parallel edges, since
$G_S$ has just three edges.
Then there are two Scharlemann cycles $\tau_1$ and $\tau_2$ of length three in $G_S$.
Let $g_i$ be the face of $G_S$ bounded by $\tau_i$ for $i=1,2$.
Clearly, $g_1$ and $g_2$ lie on opposite sides of $\widehat{T}$.
The edges of $\tau_i$ are all edges of $G_T$.
In particular, $\tau_1$ and $\tau_2$ have their edges in common.

\begin{claim} 
The edges of $\tau_i$ cannot lie in a disk in $\widehat{T}$ for $i=1,2$.
\end{claim}

\noindent
\textit{Proof of Claim 6.4.}
Suppose that the edges of $\tau_1$ (and therefore $\tau_2$) lie in a disk $D$ in $\widehat{T}$.
By the cut-and-paste operation of $g_i$, we can assume that 
$\mathrm{Int}g_i\cap D=\emptyset$ for $i=1,2$.
Then $N(D\cup V\cup g_1\cup g_2)$ gives
a connected sum of two lens spaces minus
an open $3$-ball, in $K(r)$, which is impossible.
\hfill(Proof of Claim 6.4)\qedsymbol
\bigskip

Thus there is an essential annulus $A$ in $\widehat{T}$ which contains $G_T$ by Lemma 2.1.
In particular, $G_T$ has exactly one pair of parallel edges.

\begin{claim} 
$\mathrm{Int}g_i\cap\widehat{T}=\emptyset$ for $i=1,2$.
\end{claim}

\noindent
\textit{Proof of Claim 6.5.}
Suppose that $\mathrm{Int}g_1\cap\widehat{T}\ne\emptyset$.
Let $\xi$ be an innermost component of $\mathrm{Int}g_1\cap\widehat{T}$ in $g_1$,
and let $\delta$ be the disk bounded by $\xi$ on $g_1$.
By the assumption on the loops in $S\cap T$ stated in Section 2,
$\xi$ is essential in $T$, and then
$\xi$ is parallel to $\partial A$.
We may suppose that $\delta\subset W$.
Then $\delta$ is a meridian disk of $W$.
Let $H=V\cap W$, and let $g_j\cap H\ne\emptyset$ for some $j\in\{1,2\}$.

Let $Q$ be the $2$-sphere obtained by compressing $\partial W$ along $\delta$,
and let $B$ be the $3$-ball bounded by $Q$ in $W$.
On $Q$, there is a disk $E$ which contains the edges of $\tau_j$.
After the components of $\mathrm{Int}g_j\cap Q$ are removed by the cut-and-paste 
operation of $g_j$,
$N(E\cup H\cup g_j)$ gives a punctured lens space in $B$, which is impossible.
Therefore, $\mathrm{Int}g_1\cap\widehat{T}=\emptyset$.
Similarly for $g_2$. 
\hfill(Proof of Claim 6.5)\qedsymbol
\bigskip

Now, we may assume that $g_1\subset W$ and $g_2\subset U$,
and that $H_{1,2}=V\cap W$ and $H_{2,1}=V\cap U$.
As in the proof of Lemma 5.4, let $M_1=N(A\cup H_{1,2}\cup g_1)$ and
$M_2=N(A\cup H_{2,1}\cup g_2)$.
Then $M_i$ is a solid torus, and $A$ runs three times on $M_i$ in the longitudinal
direction for $i=1,2$ by Lemma 2.1.
Let $A'_i=\mathrm{cl}(\partial M_i-A)$, then $A'_i$ is parallel to the annulus
$B=\mathrm{cl}(\widehat{T}-A)$ in $W$ if $i=1$, or $U$ if $i=2$.

\begin{claim}  
$\mathrm{cl}(U-H_{2,1})$ is a handlebody of genus two.
\end{claim}

\noindent
\textit{Proof of Claim 6.6.}
The torus $A'_2\cap B$ bounds a solid torus $U'$ in $U$, which
represents the parallelism of $A'_2$ and $B$.
Then it can be seen that $\mathrm{cl}(U-H_{2,1})$ is obtained
from $U'$ by attaching a $1$-handle $N(g_2)$.
Hence we have the desired result.
\begin{flushright}
(Proof of Claim 6.6)\qedsymbol
\end{flushright}
\bigskip

\begin{claim}  
Let $e$ be one of the parallel edges in $G_T$.
Then $e$ is an unknotting tunnel of $K$.
\end{claim}

\noindent
\textit{Proof of Claim 6.7.}
Note that $\mathrm{cl}(W-H_{1,2})$ is homeomorphic to
$\mathrm{cl}(M_1-H_{1,2})$.
Let $k=K^*\cap W$, where $K^*$ is the core of $V$.
Then $k$ is a properly embedded arc in $M_1$,
and $\mathrm{cl}(M_1-H_{1,2})=\mathrm{cl}(M_1-N(k))$.
Push $e$ into $W$ slightly.
It can be assume that $\partial e\subset\partial N(k)$.
See Figure 6.


\bigskip
$$\fbox{Figure 6}$$
\bigskip

Then it is not hard to see that
$\mathrm{cl}(M_1-N(k)\cup N(e))$ has a product structure $T\times I$.
Since $\mathrm{cl}(U-H_{2,1})$ is a handlebody of genus two by Claim 6.6,
$\mathrm{cl}(E(K)-N(e))$ is a handlebody of genus two, which gives the desired conclusion.
(Proof of Claim 6.7)\qedsymbol
\bigskip

By Lemma 6.1, $K$ is $2$-bridge, and this means that the case $r=3$ is impossible.
\end{proof}

\begin{lemma} 
The case $r=5$ is impossible.
\end{lemma}

\begin{proof}
$G_S$ has exactly five edges.
By Lemma 6.2, these edges of $G_S$
are partitioned into three families, two pairs of parallel edges and one edge which
is not parallel to the others.
However, this configuration contradicts
the fact that the faces of $G_S$ are divided into black and white sides.
\end{proof}

\bigskip
\noindent
\textit{Proof of Theorem 1.2 when $t=2$.}
By Proposition 4.5 and Theorem 5.1, the remaining cases are $r=3,5$.
These are impossible by Lemmas 6.3, 6.8.
This completes the proof of Theorem 1.2. \hfill\qedsymbol

\bigskip
\noindent
\textit{Proof of Theorem 1.3.}
Let $K$ be a genus one knot in $S^3$, and suppose that $K(r)$ is a lens space.
By Theorem 1.2, $K$ is not hyperbolic, and therefore it is either a satellite knot
or a torus knot.
If a satellite knot admits cyclic surgery, 
then it is a cable knot of a torus knot \cite{BL,Wa,Wu}.
In particular, its genus is greater than $1$. 
Thus we have that $K$ is a torus knot, and so $K$ is the trefoil.
The constraint on the slopes follows from \cite{Mo}.
The converse is obvious.\hfill\qedsymbol

\bigskip
The authors are indebted to Dr. Makoto Ozawa for helpful 
conversations. Part of this work was carried out while 
the first author was visiting at University of California, 
Davis. He would like to express hearty thanks to
Professor Abigail Thompson and the department for their hospitality.

\bibliographystyle{amsplain}

\end{document}